\documentclass{amsart}
 \usepackage{amsrefs}
  \begin{document} 
   \newtheorem{theorem}{Theorem}[section]
\newtheorem{lemma}[theorem]{Lemma}
\newtheorem{Coro}[theorem]{Corollary}
\newtheorem{Prop}[theorem]{Proposition}
\theoremstyle{definition}
\newtheorem{definition}[theorem]{Definition}
\newtheorem{example}[theorem]{Example}
\newtheorem{xca}[theorem]{Exercise}

\theoremstyle{remark}
\newtheorem{remark}[theorem]{Remark}

\numberwithin{equation}{section}

  \newcommand{\RR}{\mathbb{R}}
\newcommand{\N}{\mathbb{N}}

  \title[Maximal function for the generalized Ornstein-Uhlenbeck semigroup]
{On the maximal function for the generalized Ornstein-Uhlenbeck
semigroup.} 
\author[J. Betancor]{Jorge Betancor}
\address{Departamento de An\'alisis Matem\'atico, 
    Universidad de la Laguna,
    38271- La Laguna, Tenerife
    Islas Canarias, Spain}
 \email{jbetanco@ull.es} 
 \thanks{The first author was supported in part by Consejer\'\i a de Educaci\'on, Gobierno de Canarias, grant PI2003/068.}  
\author[L. Forzani]{Liliana Forzani}
 \address{IMAL- Facultad de Ingenier\'{\i}a Qu\'\i mica, 
    U. N. del Litoral , Guemes 3450,
    Santa Fe 3000 and CONICET, Argentina}
  \curraddr{School of Statistics, University of Minnesota, USA 
    Ford Hall 495, Minneapolis MN 55414, USA}
\email{liliana.forzani@gmail.com}
\author[R. Scotto]{Roberto Scotto}
\address{IMAL- Facultad de Ingenier\'{\i}a Qu\'\i mica, 
    U. N. del Litoral, Guemes 3450,
    Santa Fe 3000 Argentina}
\email{roberto.scotto@gmail.com}
\author[W.Urbina]{Wilfredo~O.~Urbina}
  \address{Departamento de Matem\'{a}ticas, Facultad de Ciencias,  UCV. Apt 47195,
    Los Chaguaramos, Caracas 1041-A Venezuela}
\curraddr{Department of Mathematics and Statistics, University of New Mexico,
   Albuquerque, NM 87131, USA.}
\email{wurbina@math.unm.edu}
\thanks{The fourth author was supported in part by FONACIT-Venezuela, grant G97000668.} 


\begin{abstract}
In this note we consider  the maximal function for the generalized
Ornstein-Uhlenbeck semigroup in $\RR$ associated with the
generalized Hermite polynomials $\{H_n^{\mu}\}$ and prove that it is
weak type (1,1) with respect to $d\lambda_{\mu}(x) =
|x|^{2\mu}e^{-|x|^2} dx,$ for  $\mu >-1/2$ as well as bounded on
$L^p(d\lambda_\mu) $ for $p>1$.
\end{abstract}
\maketitle

\section{Introduction and Preliminaries}

\noindent The  generalized Hermite polynomials were defined by G.
Sz\"ego in~\cite{Szego} (see problem 25, pag 380) as being orthogonal
polynomials with respect to the measure
$d\lambda(x)=d\lambda_{\mu}(x) = |x|^{2\mu}e^{-|x|^2} dx,$ with $
\mu>-1/2$. In his doctoral thesis T. S. Chihara~\cite{Chihara1955} (see
also~\cite{Chihara1978}) studied them in detail. In this paper we consider
the definition of the generalized Hermite polynonials given by M.
Rosenblum in~ \cite{Rosenblum}.

\noindent Let us denote by $H_{n}^{\mu} $ this generalized Hermite
polynomial of degree $n$, then for $n$ even
\begin{equation}\label{def3}
H_{2m}^{\mu} (x)= (-1)^m (2m) ! \frac{\Gamma(\mu +
\frac{1}{2})}{\Gamma(m+\mu + \frac{1}{2})} L^{\mu -
\frac{1}{2}}_m(x^2)
\end{equation}
and for $n$ odd
\begin{equation}\label{def4}
 H_{2m+1}^{\mu}(x)= (-1)^m (2m+1) ! \frac{\Gamma(\mu + \frac{3}{2})}{\Gamma(m+\mu + \frac{3}{2})} x L^{\mu +
 \frac{1}{2}}_m(x^2),
\end{equation}
\noindent $L^{\gamma}_m$ being the ${\gamma}$-Laguerre polynomial of
degree $m$.

\noindent Thus, for every $n\in \N$,
$$
\|H_n^\mu\|_{L^2(d\lambda)}=
\bigg(\frac{2^n(n!)^2\Gamma(\mu+1/2)}{\gamma_\mu(n)}\bigg)^{1/2},
$$
where $\gamma_{\mu}(m)$ is a generalized factorial defined by,
$$ \gamma_{\mu}(2m) = \frac{ 2^{2m} m!
\Gamma(m+\mu+\frac{1}{2})}{\Gamma(\mu+\frac{1}{2})}=(2m)!
\frac{\Gamma(m+\mu+\frac{1}{2})}{\Gamma(\mu+\frac{1}{2})}\frac{
\Gamma(\frac{1}{2})}{\Gamma(m+\frac{1}{2})},$$
$$ \gamma_{\mu}(2m+1) =\frac{ 2^{2m+1} m!
\Gamma(m+\mu+\frac{3}{2})}{\Gamma(\mu+\frac{1}{2})}=(2m)!
\frac{\Gamma(m+\mu+\frac{3}{2})}{\Gamma(\mu+\frac{1}{2})}\frac{
\Gamma(\frac{1}{2})}{\Gamma(m+\frac{3}{2})}.$$

\noindent The generalized Hermite polynomials $\{H_{n}^{\mu} \}$
have a generating function (2.5.8) of ~\cite{Rosenblum})
 which involves the  generalized
exponential function $e_{\mu}$ defined by

 \begin{equation}
e_{\mu} (z) = \sum_{m=0}^{\infty} \frac{z^m}{\gamma_{\mu}(m)}.
\end{equation}

\noindent  \noindent On the other hand each generalized Hermite
polynomial satisfies the following differential equation, see~\cite{Chihara1978},
\begin{equation}\label{difeq}
(H_{n}^{\mu})''(x) + 2(\frac{\mu}{x} - x) (H_{n}^{\mu})' (x) + 2(n -
\mu \frac{\theta_n}{x^2}) H_{n}^{\mu}(x)= 0,
 \end{equation}
 with $$ \theta_{n} =\left\{
\begin{array}
[c]{lll}%
\displaystyle{1}  & \text{if } & n \; \text{is odd}, \\
0 & \text{if } & n \; \text{is even.}%
\end{array}
\right.
$$ and $n \geq 0.$

\noindent Therefore, by considering the (differential-diference) operator
\begin{equation}\label{difop}
L_{\mu} = \frac{1}{2}\frac{d^2}{dx^2} + (\frac{\mu}{x} -
x)\frac{d}{dx}  - \mu \frac{I -\tilde{I}}{2 x^2},
\end{equation}
where $If(x)=f(x)$ and $ \tilde{I}f(x)=f(-x)$, $H_{n}^{\mu}$ turns
out to be an eigenfunction of $L_{\mu}$ with eigenvalue $-n$.

Now  we can define a Markov semigroup, see  D. Bakry~\cite{Bakry}, by
\begin{equation}\label{markvsem}
P_t(x, dy) = \sum_{n=0}^{\infty}  \frac{\gamma_{\mu}(n)}{2^n (n!)^2}
H_{n}^{\mu}(x) H_{n}^{\mu}(y) e^{-nt}\lambda(dy).
\end{equation}
This semigroup is entirely characterized
by the action on positive or bounded measurable functions by
$$T_{\mu}^t f(x)= \int_{-\infty}^{\infty} f(y) P_t(x, dy).$$
Thus the family of operators $\{T_{\mu}^t\}_{t\geq 0}$ is then a
 conservative semigroup of operators with generator $L_{\mu}$, that we will call the 
generalized Ornstein-Uhlenbeck semigroup. Therefore,
$$ \frac{ \partial T_{\mu}^t f(x)}{\partial t} = L_{\mu} T_{\mu}^t f(x).$$
 For $\mu=0$, $\{T^t_\mu\}$ reduces to the Ornstein-Uhlenbeck semigroup
whose behavior on $L^p$ was studied  by B. Muckenhoupt in~\cite{Muckenhoupt} for the
one-dimensional case. \noindent By using the generalized Mehler's formula (2.6.8) of~\cite{Rosenblum}: for $ x , \, y \in {\RR} $ and $
|z|<1 $,
\begin{equation}\label{mehler}
\sum_{n=0}^{\infty} \frac{\gamma_{\mu}(n)}{2^n (n!)^2}
H_{n}^{\mu}(x) H_{n}^{\mu}(y) z^n = \frac 1{\left( 1-z^2\right)
^{\mu+1/2}}e^{-\frac{z^2(x^2+y ^2)}{1-z^2}} e_{\mu}\left(\frac{2xy
z}{1-z^2}\right).
\end{equation} we can obtain
the following integral expression of this generalized
Ornstein-Uhlenbeck semigroup $\left\{ T^{t}_\mu\right\}$,
\begin{equation}\label{semig}
T_{\mu}^t f(x)=\frac 1{\left( 1-e^{-2t}\right) ^{\mu+1/2}}\int_{-\infty}^{\infty}e^{-\frac{%
e^{-2t}(x^2+y ^2)}{1-e^{-2t}}} e_{\mu}\left(\frac{2xy
e^{-t}}{1-e^{-2t}}\right) f(y)|y|^{2\mu}e^{-|y|^2} dy.
\end{equation}

\bigskip
In the following section we will consider the maximal operator
associated with $\{T_{\mu}^t \}_{t> 0}$, and prove it is weak type
$(1,1)$ with respect to the measure $\lambda$, bounded in  $L^\infty$
and therefore $L^p$ bounded for $1<p< \infty$ with respect to  $\lambda$. 
It is important to observe that since $\{T_{\mu}^t \}_{t> 0}$ is not a convolution
semigroup, its associated maximal operator is not bounded by the Hardy-Littlewood maximal
operator and therefore in order to prove the weak $(1,1)$ inequality with respect to
$\lambda$ it is needed to develop new techniques. The case $\mu=0$, that 
as we already said corresponds to the maximal operator of the Ornstein-Uhlenbeck semigroup,  was 
proved by Sj\"{o}gren in~\cite{Sjogren1983} in any dimension.

\noindent We will use repeatedly that
\begin{equation}\label{cuenta}
|x|^k e^{-x^2} \le C e^{-x^2/2}\le C, \ \ \forall \ x\in {\RR} .
\end{equation}
 The constant $C$ which will appear throughout this paper
may be different on each occurrence.

\medskip

\section {The maximal function of the generalized Ornstein Uhlenbeck semigroup.}
\noindent Let us define the generalized Ornstein-Uhlenbeck maximal
function as
\begin{eqnarray}
T_\mu^*f(x)& = &{\sup_{t > 0}} \;\,|T_\mu^tf(x)|,
\end{eqnarray}
 for each $x\in{\RR} $. Taking $r = e^{-t}$, we can write
$$T_\mu^*f(x)= {\sup_{0<r<1}} \left|\int_{-\infty}^{\infty}K_r(x,y)f(y)\,
d\lambda(y)\right|,$$
with
\[
K_r(x,y)=\frac{1}{\Gamma
(\mu+\frac{1}{2})(1-r^2)^{\mu+\frac{1}{2}}}e^{-(x^2+y^2)\frac{r^2}{1-r^2}}e_{\mu}(\frac{2xyr}{1-r^2}).
\]

\noindent The main result of this paper is summarized in

\begin{theorem}\label{debilpositivo}
For $\mu>-1/2$,
\begin{enumerate}
\item[i)] $T^*_\mu$ is weak type $(1,1)$ with respect to $\lambda$, i.e. there
exists a real constant $C>0$ such that for every $\eta>0$
\begin{equation}\label{desdebil}
\lambda \{ x\in {\RR} : T^*_\mu f(x)> \eta \}\le
\frac{C}{\eta }\|f\|_{1,\lambda},
\end{equation}
 where $\displaystyle
\|f\|_{1,\lambda}=\int_{{\RR} } |f(y)|d\lambda (y).$ 

\item[ii)] $T^*_\mu$ is bounded in $L^\infty $, i. e. there exists a real
constant $C>0 $ such that
\begin{equation} \label{fuerte}
|| T^*_\mu f ||_\infty \le C || f||_\infty
\end{equation}
where $\displaystyle \|f\|_\infty $ represents the $L^\infty $ norm.
\end{enumerate}
\end{theorem}

\begin{Coro}
For $\mu> -1/2$ and $p>1$,
\begin{equation}\label{desfuerte}
\|T^*_\mu f\|_{p,\lambda}\le C \,
\|f\|_{p,\lambda},
\end{equation}
 where $\displaystyle
\|f\|_{p,\lambda}^p=\int_{{\RR} } |f(y)|^p d\lambda (y).$
\end{Coro}

 \noindent This corollary follows from Marcinkiewicz interpolation theorem
 between the weak type $(1,1)$ and the
 boundedness in $L^\infty$ which will be proved in Theorem
 \ref{debilpositivo}.
 In order to prove Theorem \ref{debilpositivo} we will introduce well known bounds for the
 functions $e_\mu $ and prove two
propositions. The first one due to  I. P. Natanson and B.
Muckenhoupt (~\cite{Natanson} and~\cite{Muckenhoupt}) is a sort of a generalized
Young's inequality for Borel measures, that we will write it only
for the particular case of the measure $\lambda$ and the other one
has to do with the biggest function whose density distribution as a
function of $\eta$ with respect to $\lambda$ is bounded by $C/\eta.$

\bigskip

\noindent {\bf Properties of $e_\mu$}

\medskip

It can be proved, see (2.2.3) of~\cite{Rosenblum}, that the generalized
exponential function $e_{\mu}$ can be written as,
$$
e_\mu(x)=\Gamma(\mu+1/2)(2/x)^{\mu-1/2}(I_{\mu-1/2}(x)+I_{\mu+1/2}(x)),
$$
where $I_\nu$ denotes the modified Bessel function. Then, according
to~\cite[(2), p. 77, and (2), p. 203]{Watson}, we have the following
estimates that will be useful in the sequel
\begin{equation}
|e_\mu(x)|\le e_\mu(|x|)\le C(1+|x|)^{-\mu}e^{|x|},\,\,\,x\in
{\RR} . \label{bes}
\end{equation}
Also, $e_\mu$ admits the following integral representations
depending on the values of $\mu$ ~\cite{Rosenblum},
\begin{enumerate}
\item  if $\mu>0 \quad$ then
\begin{equation}\label{mupositivo}
e_{\mu}(x)=\frac{1}{B(\frac{1}{2},\mu)}\int_{-1}^{1}e^{xt}(1-t)^{\mu-1}(1+t)^{\mu}\,
dt,
\end{equation}
\item if $\mu=0 \quad$ then
\begin{equation}\label{mucero}
e_0(x)=e^x ,\end{equation} \item if $-\frac{1}{2}<\mu<0 \quad $ then
\begin{equation} \label{munegativo}
e_{\mu}(x)=e^x+\frac{\mu}{\mu+1/2}\frac{1}{B(1/2,\mu+1)}\int_{-1}^1(e^{xt}-e^x)(1-t)^{\mu-1}(1+t)^\mu
dt
\end{equation}
\end{enumerate}

According to (\ref{mupositivo}) it is clear that $e_\mu(x)\ge 0$,
for $ \mu \ge 0 $, $x\in {\RR} $. However, this one is not
the case when $-1/2<\mu<0$. Indeed, assume that $-1/2<\mu<0$. Since
$e^u-1\ge u$, $u>0$, we can write
\begin{eqnarray*}
e^{-x}e_\mu(x)&=&1+\frac{\mu}{\mu+1/2}\frac{1}{B(1/2,\mu+1)}\int_{-1}^1(e^{x(t-1)}-1)(1-t)^{\mu-1}(1+t)^\mu
dt\\
&\le&1-\frac{x\mu}{\mu+1/2}\frac{1}{B(1/2,\mu+1)}\int_{-1}^1(1-t)^{\mu}(1+t)^\mu
dt,\,\,\,x<0.
\end{eqnarray*}
Hence, there exists $x_0>0$ such that $e_\mu(x)<0$ for every
$x<-x_0$.

From the above we infer that the generalized Ornstein-Uhlenbeck
semigroup $\left\{ T_{\mu}^t\right\}_{t> 0}$ is a positive one when
$\mu\ge 0$ but it is not when $-1/2<\mu<0$.

\bigskip

\begin{Prop}(Natanson) \label{natprop}
Let $f$ and $g$ be two $L^1(d\lambda)$ functions. Let us assume that
g(y) is nonnegative and there is an $x\in {\RR} $ such that
$g(y)$ is monotonically increasing for $y\le x$ and monotonically
decreasing for $x\le y$, then \begin{equation}\label{natanson}
\left|\int g(y)f(y)\, d\lambda(y)\right|\leq
\|g\|_{1,\lambda}\mathcal{M}_{\lambda}f(x)
\end{equation}
\noindent where
$$\mathcal{M}_{\lambda}f(x)={\sup_{x\in
I}}\frac{1}{\lambda(I)}\int_I\left|f(y)\right|\, d\lambda(y)$$ is
the Hardy-Littlewood maximal fuction of $f$ with respect to
$\lambda$. Moreover the Hardy-Littlewood maximal fuction
$\mathcal{M}_{\lambda}f$ is weak type (1,1) and strong type (p,p)
for $p>1 $ with respect to the measure $\lambda $.
\end{Prop}
\noindent A proof of this proposition can be found in~\cite{Muckenhoupt}.

\bigskip

\begin{Prop}\label{funcionmaxima}
For $\mu>-1/2$, there is a real constant $C>0$ such that the distribution function with respect to $\lambda$ of the function
\begin{equation*}
h(x)= \max \left( \frac{1}{|x|},
|x|\right)\frac{e^{x^2}}{|x|^{2\mu}}
\end{equation*}
satisfies the inequality
$$\lambda\{x\in {\RR} : h(x)>\eta\}\le \frac{C}{\eta},$$
for any $\eta > 0$.
\end{Prop}
\begin{proof}
Since $\lambda$ is a finite measure, it is enough to prove this
result for $\eta \ge e$. Besides, due to the fact that
 $h$ is even and $\lambda$ is symmetric, then $\lambda \{ x\in
{\RR} : h(x)> \eta \}=2\lambda \{ x>0 : h(x)> \eta \}$. Now
\begin{eqnarray*}
\lambda \{ x>0: h(x)> \eta \}&\le & \lambda \left \{ 0<x<1 :
\frac{1}{x^{2\mu+1}}> \eta/e \right \} \\ & & +\lambda \left \{ x>1:
\frac{e^{x^2}}{x^{2\mu -1}}> \eta \right
\}\\&=&\int_0^{(e/\eta)^{\frac{1}{2\mu+1}}} x^{2\mu}e^{-x^2}dx\\
& &+ \int_{x_0}^{\infty} x^{2\mu}e^{-x^2}dx\\& = &I+II
\end{eqnarray*}
with $x_0>1$ and $\frac{e^{x_0^2}} {x_0^{2\mu-1}}=\eta$.  Let us
observe that $$I\le \int_0^{(e/\eta)^{1/(2\mu+1)}}
x^{2\mu}dx=\frac{e}{(1+2\mu) \eta},$$ and $$II\le C
x_0^{2\mu-1}e^{-x_0^2}= \frac{C}{\eta}.$$ For last inequality see~\cite{ForMacSco}. 
From these two bounds the conclusion of this
proposition follows.
\end{proof}

\begin{proof} \emph{of Theorem} \ref{debilpositivo}.

In order to prove this theorem it suffices to show that there exists
$C>0$ such that
\begin{equation} \lambda\{x\in
(0,\infty):T_{\mu,+}^*f(x)>\eta\}\le
\frac{C}{\eta}\|f\|_{1,\lambda},\,\,\,\eta>0,
\label{forweak}
\end{equation} and
\begin{equation}
\|T^*_{\mu, +}f\|_\infty\le C \|f\|_\infty
\label{infinity}\end{equation} for every $f\ge 0$, where
$$
T_{\mu,+}^*f(x)=\sup_{t>0}|T_{t,+}^\mu f(x)|,
$$
and
$$
T_{\mu,+}^t f(x)=\frac 1{\left( 1-e^{-2t}\right) ^{\mu+1/2}}\int_{0}^{\infty}e^{-\frac{%
e^{-2t}(x^2+y ^2)}{1-e^{-2t}}} e_{\mu}\left(\frac{2xy
e^{-t}}{1-e^{-2t}}\right) f(y)|y|^{2\mu}e^{-|y|^2} dy.
$$

Indeed, let us write $r=e^{-t}$, with $t>0$. By (\ref{bes}), we
have that
$$
K_r(x,y)\le K_r(|x|,|y|),\,\,\,x,y\in {\RR} .
$$
Then
$$
|T^t_\mu f(x)|\le T_{\mu,+}^t |f|(|x|)+T_{\mu,+}^t
|\tilde{f}|(|x|),\,\,\,x\in {\RR} ,
$$
being $\tilde{f}(x)=f(-x)$, $x\in {\RR} $. Hence,
$$
T_\mu^* f(x)\le T_{\mu,+}^* |f|(|x|)+T_{\mu,+}^*
|\tilde{f}|(|x|),\,\,\,x\in {\RR} ,
$$
and we can write, for every $\eta>0$,
\begin{eqnarray*}
\lambda\{x\in {\RR} :T_\mu^* f(x)>\eta\}&\le& \lambda\{x\in
{\RR} :T_{\mu,+}^* |f|(|x|)>\eta/2\}\\
&& \;\; +\lambda\{x\in
{\RR} :T_{\mu,+}^* |\tilde{f}|(|x|)>\eta/2\}\\
&\le&2(\lambda\{x\in
(0,\infty):T_{\mu,+}^* |f|(x)>\eta/2\}\\
&& \;\; + \lambda\{x\in (0,\infty):T_{\mu,+}^*
|\tilde{f}|(x)>\eta/2\}).
\end{eqnarray*}
Thus (\ref{desdebil}) follows from (\ref{forweak}), (\ref{infinity})
and the fact that $\|f\|_{1,\lambda}=\|\tilde{f}\|_{1,\lambda}$ and
$\|f\|_\infty= \|\tilde{f}\|_\infty$.

From now on let us assume $f\ge 0$ and $x>0$. First let us prove the
weak type $(1,1)$ inequality.
\medskip

\noindent {\bf (1)} Case $\mu=0$. This case corresponds to the
Ornstein-Uhlenbeck maximal operator which was proved to be weak type
$(1,1)$ by B. Muckenhoupt in~\cite{Muckenhoupt}.

\medskip

\noindent {\bf (2)} Case $\mu>-1/2$. By using (\ref{bes}) we can
write
\begin{eqnarray*}
T^t_{\mu,+}f(x)&\le&\frac{C}{(1-r^2)^{\mu+1/2}}\, \int_0^\infty
e^{-\frac{(x^2+y^2)r^2}{1-r^2}+\frac{2xyr}{1-r^2}}\bigg(1+\frac{2xyr}{1-r^2}\bigg)^{-\mu}f(y)\,
  d\lambda(y)\\
&=&\frac{Ce^{x^2}}{(1-r^2)^{\mu+1/2}}\, \int_0^\infty
e^{-\frac{|x-ry|^2}{1-r^2}}\bigg(1+\frac{2xyr}{1-r^2}\bigg)^{-\mu}\,
f(y)\,
d\lambda(y)\\
&=&\frac{Ce^{x^2}}{(1-r^2)^{\mu+1/2}}\,
\bigg(\int_0^{x/2r}+\int_{x/2r}^{4x/r}+\int_{4x/r}^\infty\bigg)
e^{-\frac{|x-ry|^2}{1-r^2}}\bigg(1+\frac{2xyr}{1-r^2}\bigg)^{-\mu}\,
f(y)\,
d\lambda(y)\\
&=& C(K_{1,r}f(x)+K_{2,r}f(x)+K_{3,r}f(x)).
\end{eqnarray*}

Let us observe that if $0<y<x/2r$, then $x-ry>x/2$ and
$$\frac{1}{(1-r^2)^{\mu+1/2}}\left(1+\frac{2rxy}{1-r^2}\right)^{-\mu}\leq
\frac{1}{(1-r^2)^{\mu+1/2}}+\frac{x^{-2\mu}}{(1-r^2)^{1/2}}, $$ thus
$$
K_{1,r}f(x)\le C e^{x^2}
\left(\frac{1}{(1-r^2)^{\mu+1/2}}+\frac{x^{-2\mu}}{(1-r^2)^{1/2}}\right)e^{-\frac{x^2}{4(1-r^2)}}\|f\|_{1,\lambda}\le
C\frac{e^{x^2}}{x^{2\mu+1}}\|f\|_{1,\lambda},
$$
where last inequality is obtained as an application of
(\ref{cuenta}).

 On the other hand, if $y>\frac{4x}{r}$, then
$ry-x>x$, and again by applying (\ref{cuenta}) repeatedly in the
sequel below
\begin{eqnarray*}
\frac{e^{-\frac{|x-ry|^2}{1-r^2}}}{(1-r^2)^{\mu+1/2}}\bigg(1+\frac{2rxy}{1-r^2}\bigg)^{-\mu}&=&
\frac{e^{-\frac{|x-ry|^2}{1-r^2}}}{(1-r^2)^{\mu+1/2}}\bigg(1+\frac{2x
(ry-x)+2x^2}{1-r^2}\bigg)^{-\mu}\\&\le & C \
e^{-\frac{x^2}{2(1-r^2)}}
\left(\frac{1}{(1-r^2)^{\mu+1/2}}+\frac{x^{-\mu}}{(1-r^2)^{\frac{\mu+1}{2}}}+
\right. \\ & & \hspace{2.5cm} \left.
\frac{x^{-2\mu}}{(1-r^2)^{1/2}}\right)\\ &\le &
\frac{C}{x^{2\mu+1}},
\end{eqnarray*}
we get
$$
K_{3,r}f(x)\le C\frac{e^{x^2}}{x^{2\mu+1}}\|f\|_{1,\lambda}.
$$
Finally for $\frac{x}{2r}\le y\le \frac{4x}{r}$ we have the
following estimate
\begin{equation}\label{ultima}
\frac{1}{(1-r^2)^{\mu+1/2}} \left(1+\frac{2rxy}{1-r^2}\right)^{-\mu}
\leq \frac{1}{x^{2\mu+1}}+\frac{x^{-2\mu}}{(1-r^2)^{1/2}},
\end{equation}
which is immediate for $\mu\ge 0$ and for $\mu<0$ one has to argue
between $\frac{2rxy}{1-r^2}\le 1$ and its complement. Now by taking
into account inequality (\ref{ultima}) we are ready to estimate
$K_{2,r}f(x)$ and for that we consider two cases. If $0<r\le 1/2$ we
have
$$
K_{2,r}f(x)\le C
\left(\frac{1}{x}+1\right)\frac{e^{x^2}}{x^{2\mu}}\|f\|_{1,\lambda},
$$
and, if $1/2<r<1$ then
$$
K_{2,r}f(x)\le
C\left(\frac{e^{x^2}}{x^{2\mu+1}}\|f\|_{1,\lambda}+\frac{e^{x^2}}{(1-r^2)^{1/2}x^{2\mu}}\int_0^\infty
 N(r,x,y)f(y)d\lambda(x)\right),
$$
with
\begin{equation}\label{nucleonatanson} N(r,x,y)=\left\{
\begin{array}{lcl}
1& \mbox{if} & y\in \left[x, \frac{x}{r}\right]\\
e^{-\frac{|x-ry|^2}{1-r^2}}& \mbox{if} & y\in \left[\frac{x}{2r},
\frac{4x}{r}\right]\setminus \left[x, \frac{x}{r}\right]\\
0& &\mbox{otherwise}.
\end{array}\right.
\end{equation}  Since $N(r,x,.)$ is a Natanson kernel (see (\ref{natanson})), we get
$$K_{2,r} f(x)\le
C\left(\frac{e^{x^2}}{x^{2\mu+1}}\|f\|_{1,\lambda}+
\frac{e^{x^2}}{x^{2\mu}(1-r^2)^{1/2}}\|N(r,x,.)\|_{1,\lambda}\,
\mathcal{M}_{\lambda}f(x)\right).$$ Let us prove that
\begin{equation}\label{acotacion}\|N(r,x,.)\|_{1,\lambda}\le C x^{2\mu}
(1-r^2)^{1/2}e^{-x^2}.\end{equation} Indeed,
\begin{eqnarray*}
\int_{{\RR} }N(r,x,y)\,
d\lambda(y)&=&\int_x^{x/r}e^{-y^2}y^{2\mu}\, dy +
\int_{x/2r}^{x}e^{-\frac{|x-ry|^2}{1-r^2}} e^{-y^2}y^{2\mu}\, dy
\\ & & +\int_{x/r}^{4x/r}e^{-\frac{|x-ry|^2}{1-r^2}} e^{-y^2}y^{2\mu}\,
dy\\ &\sim& x^{2\mu}\left( \int_x^{x/r}e^{-y^2}\, dy + e^{-x^2}
\int_{x/2r}^{x}e^{-\frac{|rx-y|^2}{1-r^2}}\, dy \right.\\
& &\left. + \; e^{-x^2}\int_{x/r}^{4x/r}e^{-\frac{|rx-y|^2}{1-r^2}}
\, dy\right)
\\ & \le & C\, x^{2\mu} e^{-x^2}\left(\min \left(\frac{1}{x},
(1-r)x\right) \right.\\ & &
\left. + \; \int_{{\RR} }e^{-\frac{|rx-y|^2}{1-r^2}}\, dy\right)\\
&\le& C\, x^{2\mu} (1-r^2)^{1/2}e^{-x^2}.
\end{eqnarray*}
 Now gathering together all the bounds
obtained above, we get
$$
T_{\mu,+}^t f(x)\le C (h(x) \|f\|_{1,\lambda}+ {\mathcal M}_\lambda
f (x)),
$$
for all $t>0$, where $h$ is the function defined in Proposition
\ref{funcionmaxima}. Thus the weak type $(1,1)$ of $T^{*}_{\mu,+}$
follows from propositions \ref{natprop} and \ref{funcionmaxima}.

\medskip

\noindent Now let us take care of the boundedness of $T^*_{\mu,+}$
in $L^\infty.$

\noindent For the case $\mu\ge 0$ this boundedness is immediate
since its kernel is non-negative and its integral equals 1.
Therefore let us study just the case $-1/2<\mu <0.$ By using
(\ref{bes}) and proceeding like in case 2 of the weak type $(1,1)$
inequality
\begin{eqnarray*}
T^*_{\mu, +}f(x)&\le & \frac{C}{(1-r^2)^{\mu+1/2}}\int_0^\infty
e^{-\frac{(x^2+y^2)r^2}{1-r^2}+\frac{2xyr}{1-r^2}}\bigg(1+\frac{2xyr}{1-r^2}\bigg)^{-\mu}f(y)\,
  d\lambda(y)\\&\le&\frac{C}{(1-r^2)^{\mu+1/2}}\, \int_0^\infty
e^{-\frac{|rx-y|^2}{1-r^2}}\bigg(1+\frac{2xyr}{1-r^2}\bigg)^{-\mu}\,
y^{2\mu}\, dy \  \|f\|_\infty\\ &=& \frac{C}{(1-r^2)^{\mu+1/2}}\,
\int_0^\infty
e^{-\frac{|rx-y|^2}{1-r^2}}\bigg(1+\frac{2(rx-y)y}{1-r^2}+\frac{2y^2}{1-r^2}\bigg)^{-\mu}\,
y^{2\mu}\, dy \  \|f\|_\infty \\ &\le& C \bigg( \int_0^\infty
\frac{e^{-\frac{|rx-y|^2}{1-r^2}}}{(1-r^2)^{\mu+1/2}}\bigg(1+\frac{2|rx-y|y}{1-r^2}\bigg)^{-\mu}\,
y^{2\mu}\, dy \\ & &+ \int_0^\infty
\frac{e^{-\frac{|rx-y|^2}{1-r^2}}}{(1-r^2)^{1/2}}dy\bigg)\|f\|_\infty
\end{eqnarray*}
In order to prove that the first integral of last inequality is
bounded by a constant independent of $r$, $y$, and $x$ first
 we use (\ref{cuenta}) to get the inequality
$$\bigg(\frac{2|rx-y|y}{1-r^2}\bigg)^{-\mu}e^{-\frac{|rx-y|^2}{1-r^2}}\le
C
\bigg(\frac{y}{(1-r^2)^{1/2}}\bigg)^{-\mu}e^{-\frac{|rx-y|^2}{2(1-r^2)}},$$
then we split the integral in two subintervals one from $0$ to
$\sqrt{1-r^2}$ and the other from $\sqrt{1-r^2}$ to $\infty$  and we
call them $I$ and $II$. Now we proceed to bound each part.
\begin{eqnarray*}
I& = & \int_0^{\sqrt{1-r^2}}
\frac{e^{-\frac{|rx-y|^2}{2(1-r^2)}}}{(1-r^2)^{\mu+1/2}}\bigg(1+\bigg(\frac{y}{(1-r^2)^{1/2}}
\bigg)^{-\mu}\bigg)\, y^{2\mu}\, dy \\&\le&
\int_0^{\sqrt{1-r^2}}\frac{y^{2\mu}}{(1-r^2)^{\mu+1/2}}\,
dy+\int^{\sqrt{1-r^2}}_0 \frac{y^{\mu}}{(1-r^2)^{(\mu+1)/2}}\, dy\le
C,
\end{eqnarray*}
and
\begin{eqnarray*}
II& = & \int^\infty_{\sqrt{1-r^2}}
\frac{e^{-\frac{|rx-y|^2}{2(1-r^2)}}}{(1-r^2)^{\mu+1/2}}\bigg(1+\bigg(\frac{y}{(1-r^2)^{1/2}}
\bigg)^{-\mu}\bigg)\, y^{2\mu}\, dy \\&\le&
\int^\infty_{\sqrt{1-r^2}}\frac{e^{-\frac{|rx-y|^2}{2(1-r^2)}}}{(1-r^2)^{\mu+1/2}}
(\sqrt{1-r^2})^{2\mu}\, dy+\int_{\sqrt{1-r^2}}^\infty
\frac{e^{-\frac{|rx-y|^2}{2(1-r^2)}}}{(1-r^2)^{\mu+1/2}}
\frac{y^{\mu}}{(1-r^2)^{-\mu/2}}\, dy \\ &\le& 2 \int_0^\infty
\frac{e^{-\frac{|rx-y|^2}{2(1-r^2)}}}{(1-r^2)^{1/2}}dy \le C.
\end{eqnarray*}
This ends the proof of the boundedness of $T^*_{\mu,+}$ in
$L^\infty$ and at the same time the proof of Theorem
\ref{debilpositivo}.

\end{proof}

\end{document}